\documentclass[11pt]{amsart}
\usepackage{a4,amscd,amssymb,amsfonts,verbatim,times,graphicx,pstricks}

\newcommand{\R}{{\mathbb R}}
\newcommand{\Za}{{\mathbb Z}}
\renewcommand{\phi}{\varphi}
\newcommand{\half}{\frac{1}{2}}
\renewcommand{\Re}{\operatorname{Re}}

\newcommand{\tM}{\widetilde{M}}
\newcommand{\tg}{\widetilde{g}}
\newcommand{\spec}{\operatorname{spec}}
\newcommand{\mmin}{\mathcal{M}_{\rm{min}}}
\newcommand{\Ahat}{\widehat{A}}
\newcommand{\Mhat}{\widehat{M}}
\newcommand{\scal}{{\rm scal}}
\newcommand{\dvol}{{\rm dvol}}
\newcommand{\ind}{\operatorname{ind}}

\newtheorem{thm}{Theorem}[section]
\newtheorem{lemma}[thm]{Lemma}
\newtheorem{prop}[thm]{Proposition}
\newtheorem{cor}[thm]{Corollary}

\newtheorem{remark}[thm]{Remark}
\newtheorem{definition}[thm]{Definition}
\newtheorem{notation}[thm]{Notation}
\newcommand{\Remark}[1]{\begin{remark}{\rm #1}\end{remark}}
\newcommand{\Definition}[1]{\begin{definition}{\rm #1}\end{definition}}

\newcommand{\f}[1]{$ \, #1  \, $}
\newcommand{\geese}[1]{``{#1}''}                                               
\begin{document}

\title[Surgery and the Spectrum of the Dirac Operator]
{Surgery and the Spectrum of the Dirac Operator}

\author{Christian B\"ar and Mattias Dahl}

\address{Universit\"at Hamburg\\
FB Mathematik\\
Bundesstr.~55\\
20146 Hamburg\\
Germany}

\email{baer@math.uni-hamburg.de, mattias.dahl@math.uni-hamburg.de}


\date{\today}

\begin{abstract}
We show that for generic Riemannian metrics on a simply-connected closed 
spin manifold of dimension $\ge 5$ the dimension of the space of harmonic 
spinors is no larger than it must be by the index theorem. The same result 
holds for periodic fundamental groups of odd order. 

The proof is based on a surgery theorem for the Dirac spectrum which says 
that if one performs surgery of codimension $\ge 3$ on a closed Riemannian 
spin manifold, then the Dirac spectrum changes arbitrarily little provided 
the metric on the manifold after surgery is chosen properly.
\end{abstract}

\maketitle

\setcounter{section}{-1}

\section{Introduction} 
Classical Hodge-deRham theory establishes a tight link between the 
analysis of the Laplace operator acting on differential forms of a 
compact Riemannian manifold and its topology. Specifically, the 
dimension of the space of harmonic $k$-forms is a topological invariant, 
the $k^{\mbox{\footnotesize th}}$ Betti number.

The question arises whether a similar relation holds for other elliptic 
geometric differential operators such as the Dirac operator on a compact 
Riemannian spin manifold. It is not hard to see that the dimension $h_g$ 
of the space of harmonic spinors is a conformal invariant, it does not 
change when one replaces the Riemannian metric $g$ by a conformally 
equivalent one \cite[Prop.~1.3]{hitchin74a}. Moreover, the Atiyah-Singer 
index theorem implies a topological lower bound on $h_g$.

Berger metrics on spheres of dimension $4k+3$ provide examples showing 
that in general $h_g$ depends on the metric and is not topological, see 
\cite[Prop.~3.2]{hitchin74a} and \cite[Thm.~3.1]{baer96b}.
Also for surfaces  of genus at least 3 the number $h_g$ varies with the 
choice of metric \cite[Thm.~2.6]{hitchin74a}. 
All known examples indicate that the following two conjectures should be true. 
On the one hand, we should have

{\bf Conjecture A.}
{\em Harmonic spinors are not topologically obstructed, i.e.\ on any compact 
spin manifold of dimension at least three there is a metric $g$ 
such that $h_g>0$.}

Perhaps $h_g$ is even unbounded. Conjecture A has been shown to be true 
in dimensions $n \equiv 0,1,3 \mod 8$ by Hitchin \cite[Thm.~4.5]{hitchin74a} 
using a topological approach, while the first author \cite[Thm.~A]{baer96b} 
proves it for $n \equiv 3,7 \mod 8$ using an analytic approach. 
Conjecture A has been shown to hold for 
spheres of dimension $n \equiv 0 \mod 4$ with unbounded $h_g$ by Seeger 
\cite{seeger00a}.

Conjecture A is not true in dimension 2.
For example, on $S^2$ all metrics are conformally equivalent, 
hence $h_g$ cannot change with the metric.
In fact, $h_g = 0$ and one even has a simple geometric lower bound for all
Dirac eigenvalues $\lambda$ on the 2-sphere equipped with any metric
$
\lambda^2 \ge 4\pi/\text{area}(M),
$
see \cite[Thm.~2]{baer92d}.

In this paper we do not deal with the rather exceptional metrics
giving rise to large $h_g$ but instead we study generic metrics. Versions 
of the second conjecture can be found in  
\cite{baer96b,bourguignon-gauduchon92a,kotschick96a}.

{\bf Conjecture B.}
{\em On any compact connected spin manifold for a generic metric $h_g$ is no 
larger than it is forced to be by the index theorem.}

This means that up to a ``small'' set of exceptional metrics an analogue 
of classical Hodge-deRham theory also holds for the Dirac operator. 
In particular, in dimensions divisible by 4 Conjecture B claims that for
a generic metric either the kernel or the cokernel of the chiral Dirac
operator vanishes.
Using a variational approach Maier \cite{maier97a} proved Conjecture~B for 
dimension $n \le 4$. In the present article we deal with the case of 
dimension $n\ge 5$. We show

{\bf Theorem~\ref{thm-simplyconnectedinmmin}.}
{\em Conjecture~B is true for all simply connected spin 
manifolds of dimension at least five.}

Our approach, very different from Maier's, follows the strategy that has 
successfully been used to characterize manifolds admitting metrics of 
positive scalar curvature. We show that the class of manifolds for which 
Conjecture~B holds is closed under surgery of codimension $\ge 3$. Then 
we find a set of manifolds for which the conjecture holds and which
generates the spin bordism ring $\Omega_\ast^{\text{\footnotesize spin}}$.
This set consists of manifolds of positive scalar curvature and a few more
special manifolds. 
Results from bordism theory then imply the theorem.

Using more refined spin bordism rings one can also deal with certain 
classes of non-simply connected manifolds. As an example, we show

{\bf Theorem~\ref{thm-periodic}.}
{\em Suppose $M$ is a compact connected spin manifold with $\dim M \geq 5$
and fundamental group $\pi$ a periodic group of odd order. 
Then Conjecture~B holds for $M$.}

By exhibiting concrete examples Kotschick \cite{kotschick96a} shows that 
the analogue of Conjecture~B does not hold in the class of K\"ahler 
manifolds.

The most involved part of our technique is to show the invariance under
surgery mentioned above. To get this we show that if one performs surgery 
of codimension $\ge 3$ to a closed Riemannian spin manifold, then the 
lower part of the Dirac spectrum changes arbitrarily little provided the 
metric on the manifold after surgery is chosen properly. More precisely, 
we show

{\bf Theorem~\ref{thm-spectralsurgery}.}
{\em Let $\, (M, g) \,$ be a closed Riemannian manifold
equipped with a spin structure. Let $\, N \subset M \,$ be an embedded 
sphere of codimension $\, k \ge 3 \,$ and with trivialized tubular 
neighborhood. Let $\, \widetilde{M} \,$ be obtained from \f{M} by surgery 
along \f{N} together with the resulting spin structure. Let 
$\, \varepsilon > 0 \,$ and 
$\, \Lambda > 0 \,,\; \pm \Lambda \not\in {\rm spec} \, (D_g) \,$.
Then there exists a Riemannian metric $\, \widetilde{g} \,$ on 
$\, \widetilde{M} \,$ such that $\, D_g \,$ and $\, D_{\widetilde{g}} \,$
are $\, (\Lambda, \varepsilon)$-spectral close.}

Here $(\Lambda, \varepsilon)$-spectral close means that the Dirac 
eigenvalues of $M$ and $\, \widetilde{M} \,$ in the interval $(-\Lambda,
\Lambda)$ differ at most by $\varepsilon$. 

{\bf Acknowledgement.}
We would like to thank J.~Rosenberg and S.~Stolz for valuable hints.
Many thanks also go to the referee for helpful suggestions on how to
improve the presentation of the material.

The first author was partially supported by the Research and Training 
Network ``EDGE''. The second author was supported by a Post-doctoral 
scholarship from the Swedish Foundation for International Cooperation 
in Research and Higher Education ``STINT''. Both authors have also been 
partially supported by the Research and Training Network 
``Geometric Analysis'' funded by the European Commission.

\section{Surgery and the dirac operator}
By a {\em Riemannian spin manifold} we mean a Riemannian manifold 
$(M,g)$ together with the choice of a spin structure. Hence we can 
form the {\em spinor bundle} $\Sigma M$ which is a Hermitian vector bundle
over $M$ of rank $2^{[n/2]}$ where $n$ is the dimension of $M$.
For even $n$ the spinor bundle splits as 
\begin{equation}
\Sigma M = \Sigma^+ M \oplus \Sigma^-M
\label{chiralitysplitting}
\end{equation}
into two subbundles of equal rank. They are called 
{\em half-spinor bundles} of positive and negative {\em chirality}.
Sections of $\Sigma M$ are {\em spinors}. The classical 
{\em Dirac operator}, sometimes also called 
{\em Atiyah-Singer operator}, is a formally self-adjoint elliptic 
differential operator of first order acting on the space of spinors.
We denote the Dirac operator by $D$ and if we want to emphasize the
dependence on the Riemannian metric we write $D_g$. In even dimensions 
it interchanges chiralities, that is with respect to splitting 
(\ref{chiralitysplitting}) it is of the form
$$
D_g = 
\begin{pmatrix}
0  & D_g^- \\ 
D_g^+ & 0
\end{pmatrix}.
$$
For any differentiable function $f$ and spinor $\varphi$ on $M$ we 
have the formula
$$
D(f\varphi) = fD\varphi + \nabla f \cdot \varphi
$$
where $\nabla f$ is the gradient of $f$ and $X\cdot \varphi$ denotes
{\em Clifford multiplication} of a tangent vector $X$ with a spinor 
$\varphi$. Another important formula concerns the square of the Dirac 
operator,
$$
D^2 = \nabla^\ast\nabla + \frac{1}{4} \, {\rm scal}
$$
where $\nabla$ is the Hermitian connection on $\Sigma M$ induced by 
the Levi-Civita connection and scal denotes scalar curvature. This is 
usually called the {\em Lichnerowicz formula} \cite{lichnerowicz63a}
but was already known to Schr\"odinger in the 1930's 
\cite{schroedinger32a}.

We will always consider closed manifolds $M$. Standard elliptic theory 
then tells us that the spectrum of $D$ is real and discrete. 
In particular, $D_g^+$ is a Fredholm operator and its index is given by the
celebrated Atiyah-Singer index formula
$$
\ind(D_g^+) = \Ahat(M)
$$
where $\Ahat(M)$ is the {\em $\Ahat$-genus} of $M$ constructed out of 
the Pontrjagin classes of $M$. There is a generalization of this index 
theorem to which we will return in Section~\ref{sec-harmonic}.
See \cite{lawson-michelsohn89a} for a thorough introduction to spin geometry.

Next let us recall the concept of {\em surgery}.
Let $N$ be an embedded sphere of codimension $k$ in $M$ with a 
trivialized neighborhood, meaning there is a diffeomorphism $\Phi$ 
mapping $S^{n-k} \times D^k$ diffeomorphically onto its image in $M$ 
such that $S^{n-k} \times \{0\}$ is mapped onto $N$. Surgery now 
consists of removing $\Phi(S^{n-k} \times D^k)$ and gluing in 
$D^{n-k+1}\times S^{k-1}$ along the common boundary 
$S^{n-k}\times S^{k-1}$. We will only perform surgeries in codimension 
$k\ge3$. If $n-k \not= 1$, then $S^{n-k}\times S^{k-1}$ is simply 
connected (for $n-k=0$ the two components are simply connected) and carries a 
unique spin structure. This spin structure 
extends to the unique spin structure on $D^{n-k+1}\times S^{k-1}$.
Thus any spin structure on $M$ induces one on the boundary of 
$\Phi(S^{n-k} \times D^k)$ which then extends uniquely to one on 
$\tM$, the manifold after surgery. If $n-k=1$, then 
$S^{n-k}\times S^{k-1}$ has two different spin structures but only 
one of them extends to $D^{n-k+1}\times S^{k-1}$. In this case it 
will be assumed that the trivialization $\Phi$ of the neighborhood 
of $N$ is chosen such that it induces on the boundary the spin structure 
which extends. Then again we obtain a spin structure on $\tM$.

The main theorem will roughly say that the spectrum of the Dirac operator
in an arbitrarily fixed range will change only very little if one performs
surgery in codimension $\ge 3$ and chooses the metric on the resulting
manifold $\tM$ carefully. 
To formulate this in a precise way we make the following 
definition.

\Definition{Let $\, \varepsilon > 0 \,$ and $\, \Lambda > 0 \,$.
We call two operators with discrete spectrum {\em $\, (\Lambda, 
\varepsilon)$-spectral close} if 
\begin{itemize}
\item
$\pm\Lambda$ are not eigenvalues of either operator.
\item
Both operators have the same total number $m$ of eigenvalues in the 
interval $(-\Lambda, \Lambda)$.
\item
If the eigenvalues in $(-\Lambda, \Lambda)$ are denoted by 
$\lambda_1 \le \cdots \le \lambda_m$ and $\mu_1 \le \cdots \le \mu_m$ 
respectively (each eigenvalue being repeated according to its 
multiplicity), then $|\lambda_j-\mu_j|<\varepsilon$ for $j=1,\ldots,m$.
\end{itemize}
}

\begin{thm} \label{thm-spectralsurgery}
Let $\, (M, g) \,$ be a closed Riemannian manifold
equipped with a spin structure. Let $\, N \subset M \,$ be an embedded 
sphere of codimension $\, k \ge 3 \,$ and with trivialized tubular 
neighborhood. Let $\, \widetilde{M} \,$ be obtained from \f{M} by surgery 
along \f{N} together with the resulting spin structure. Let 
$\, \varepsilon > 0 \,$ and 
$\, \Lambda > 0 \,,\; \pm \Lambda \not\in {\rm spec} \, (D_g) \,$.
Then there exists a Riemannian metric $\, \widetilde{g} \,$ on 
$\, \widetilde{M} \,$ such that $\, D_g \,$ and $\, D_{\widetilde{g}} \,$
are $\, (\Lambda, \varepsilon)$-spectral close.
\end{thm}

\Remark{
If \f{U} is an open neighborhood of \f{N} in \f{M}, then 
$\, \widetilde{M} \,$ is of the form
$$
\widetilde{M} \,=\, (M \setminus U) \cup \widetilde{U}
$$
where $\, \widetilde{U} \,$ is an open subset of $\, \widetilde{M} \,$.
The proof will show that for any such \f{U} the metric 
$\, \widetilde{g} \,$ can be chosen in such a way that it coincides 
with \f{g}on $\, M \setminus U \,$.
}

By choosing $\Lambda_i>0$ such that $\Lambda_i \to \infty$ and 
$\varepsilon_i>0$ such that $\varepsilon_i \to 0$ for $i \to \infty$
we obtain the following corollary.

\begin{cor}
Let $\, (M, g) \,$ be a closed Riemannian manifold
equipped with a spin structure. Let $\, N \subset M \,$ be an embedded 
sphere of codimension $\, k \ge 3 \,$ and with trivialized tubular 
neighborhood. Let $\, \widetilde{M} \,$ be obtained from \f{M} by surgery 
along \f{N} together with the resulting spin structure.
Then there is a sequence of Riemannian metrics $\widetilde{g}_i$ on 
$\widetilde{M}$ such that the Dirac eigenvalues of $(\widetilde{M},
\widetilde{g}_i)$ converge to exactly the Dirac eigenvalues of $(M,g)$.
\end{cor}

Again we note that the metrics $\widetilde{g}_i$ can be chosen such that
they converge to $g$ away from the sphere $N$ in the following sense:
To any compact subset $K\subset M$ which does not intersect $N$ there
is an isometric copy of $K$ contained in $(\widetilde{M},\widetilde{g}_i)$
for $i$ sufficiently large.

The next section will be devoted to the proof of 
Theorem~\ref{thm-spectralsurgery}. Since the proof is somewhat involved 
the reader who is primarily interested in applications of the theorem 
may continue directly with Section~\ref{sec-harmonic} at a first reading.

\section{Proof of the surgery theorem}
The proof of Theorem~\ref{thm-spectralsurgery} consists of two
parts.
In the first part, which is quite simple, we show that to any Dirac 
eigenvalue of the original manifold $M$ in the range $(-\Lambda,\Lambda)$ 
there is a nearby eigenvalue of the manifold $\tM$ after surgery. This is 
done by multiplying a corresponding eigenspinor with a suitable cut--off 
function and plugging the resulting spinor into the Rayleigh quotient for 
the Dirac operator on $\widetilde{M}$. Here the choice of metric on 
$\widetilde{M}$ is not so crucial. It suffices that the metric coincides 
with the original metric outside a sufficiently small neighborhood of the 
sphere along which the surgery is performed.

The second, more complicated, part consists of showing that 
$\widetilde{M}$ has no further Dirac eigenvalues in the interval 
$(-\Lambda,\Lambda)$. This requires a finer analysis of the behavior 
of eigenspinors. 
First we make the scalar curvature on $M$ very large on a neighborhood $U$ 
of the surgery sphere by a $C^1$-small deformation of the metric 
(Section 2.1). 
This property persists under surgery of codimension $\ge 3$. 
Here we use the construction of 
Gromov and Lawson which shows that the class of manifolds 
admitting metrics of positive scalar curvature is closed under surgery of 
codimension $\ge 3$. Hence scalar curvature is very large in the 
corresponding  open subset $\widetilde{U}$ of $\tM$. This implies that 
most of the $L^2$-norm of an eigenspinor is supported in the complement 
$\tM \setminus \widetilde{U} = M \setminus U$ (Section 2.2). 
Then we do a finer analysis of the distribution of the $L^2$-norm on
certain annular regions (Section 2.3). 
This finally allows us to multiply an eigenspinor 
of $\widetilde{M}$ with a cut-off function whose gradient is supported 
in such an annular region where the eigenspinor carries only very little 
$L^2$-norm.
Therefore the bad error term in the Rayleigh quotient coming from the
gradient of the cut-off function is under good control. 
It follows that to each eigenvalue of $\widetilde{M}$ there
corresponds one of $M$ and we are done.

It should be noted that the proof is not at all symmetric in $M$ and 
$\widetilde{M}$ because the metric on $\widetilde{M}$ is not given 
beforehand and therefore various geometric quantities are not under 
a-priori control. To make the second part work we must choose the metric 
on $\widetilde{M}$ very carefully.

	\subsection{Increasing scalar curvature along submanifolds}

In this subsection we show that the scalar curvature can be made large near
a submanifold by a $C^1$-small deformation of the metric (which does not
change the Dirac eigenvalues much).
Then Lemma~\ref{distributionlemma} of the next subsection will show that 
only little $L^2$-norm of eigenspinors is contained in the 
problematical region near the sphere along which we will perform surgery.

\begin{prop}
\label{prop-incrscalcurv}
Let $\, (M, g) \,$ be an $n$-dimensional Riemannian manifold, let
$\, N \subset M \,$ be a compact submanifold of positive codimension.
Let \f{U} be a compact neighborhood of \f{N} in \f{M}. Then there exists 
a sequence of smooth Riemannian metrics $\, g_j \,$ on \f{M} such that
\begin{itemize}
\item 	$\, g_j \,$ is conformally equivalent to \f{g},
\item 	$\, g_j = g \,$ on $\, M \setminus U \,$,
\item 	$\, g_j \rightarrow g \,$ in the $C^1$-topology as
	$\, j \to \infty \,$,
\item 	$\, \min\limits_{N} \, {\rm scal}_{ g_j} \to \infty \,$ as
	$\, j \to \infty \,$,
\item	there exists a constant $\, S_0 \in \R \,$ such that on \f{U} 
	and for all \f{j}
	$$
	{\rm scal}_{ g_j} \ge S_0 \;.
	$$
\end{itemize}
\end{prop}

\begin{proof}
The scalar curvature of a conformally equivalent metric
$\widehat{g} \,=\, e^{2 u} \, g$ is given by
$$
{\rm scal}_{ \widehat{g}} \,=\, e^{- 2 u} \,
\left(\, {\rm scal}_{g} + 2 (n - 1) \, \Delta u 
- (n - 2) \, (n - 1) | d u |^2 \,\right) \;,
$$
see \cite[Thm.~1.159]{besse87a}. Hence it is sufficient to find a 
sequence of smooth functions
$\, u_j \,$ on \f{M} such that
\begin{itemize}
\item	$\, u_j = 0 \,$ on $\, M \setminus U \,$,
\item	$\, u_j \to 0 \,$ uniformly,
\item	$\, d u_j \to 0 \,$ uniformly,
\item	$\, \Delta u_j \,$ uniformly bounded from below and
\item	$\, \min\limits_{N} \,\Delta u_j \to \infty \,$ as
	$\, j \to \infty \,$.
\end{itemize}

To construct such functions we first assume that the normal bundle
of \f{N} in \f{M} contains a trivial line bundle, that is we assume 
there is a smooth unit normal field \f{\nu} on \f{N}. Choose a relatively 
compact neighborhood $\, U_1 \,$ of \f{N} in \f{U} sufficiently small so 
that the Riemannian exponential map maps a neighborhood 
$\, \widehat{U} \,$ of the zero section in the normal bundle 
diffeomorphically onto $\, U_1 \,$,
$$
\exp \, |_{\widehat{U}} : \widehat{U} \stackrel{\approx}{\longrightarrow}
U_1 \subset\subset U \subset M \;.
$$
We define a smooth function \f{\rho} by $\rho \, (x) \,=
\, g \, ( (\left. \exp \right|_{\widehat{U}})^{- 1} \, (x), \nu) \;$ on 
$\, U_1 \,$.
On $\, M \setminus U \,$ we define \f{\rho} to be $1$, and we extend 
\f{\rho} in some way to $\, U \setminus U_1 \,$ so that we obtain a 
smooth function defined on all of \f{M} with the following properties.
\begin{itemize}
\item	$\, \rho \equiv 1 \,$ on $\, M \setminus U \,$,
\item	$\, \rho \equiv 0 \,$ on \f{N},
\item	$\, {\rm grad} \, \rho = \nu \,$ on \f{N}, in particular
	$\, | d \rho | \equiv 1 \,$ on \f{N},
\item	$\, | d \rho | \le C_1 \,$ on \f{M},
\item	$\, | \Delta \rho | \le C_2 \,$ on \f{M}.
\end{itemize}

Next, we pick continuous functions 
$\, \widehat{\varphi}_j : \R \rightarrow \R \,$,
smooth on $\, \R  \setminus \{ 0 \} \,$ such that
\begin{itemize}
\item	$\widehat{\varphi}_j \, (0) =  \frac{1}{j}$,
\item	$0 \le \widehat{\varphi}_j \le  \frac{1}{j}$ on $\R$,
\item	$\widehat{\varphi}_j \equiv 0 \,$ on $\R \setminus (- 1, 1)$,
\item	$\widehat{\varphi}_{j}^{\prime} = \frac{2}{j}$ on 
	$ [\, -  \frac{1}{4}, 0) $,  
	$\widehat{\varphi}_{j}^{\prime} = - \frac{2}{j}$ on 
	$ (0,  \frac{1}{4} \,]$,
\item	$| \widehat{\varphi}^{\prime}_{j} | \le  \frac{2}{j}$ on 
	$ \R \setminus \{ 0 \} $,
\item	$| \widehat{\varphi}_{j}^{\prime\prime} | \le  \frac{C_3}{j}$
        on $\R\setminus \{0\}$. 
\end{itemize}

\begin{center}
\includegraphics[width=12cm]{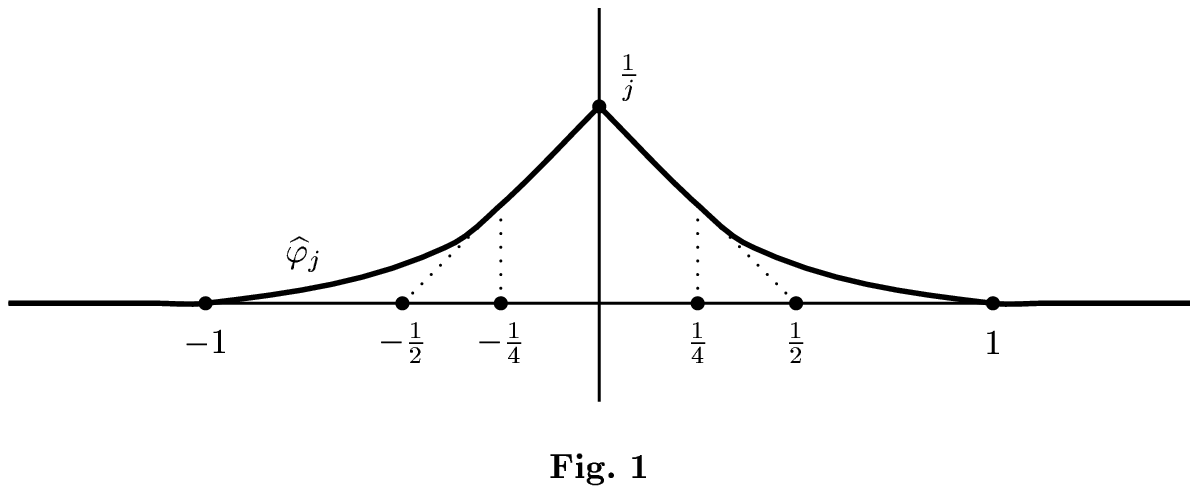}
\end{center}

Note that $\, \widehat{\varphi}_j \,$ can simply be chosen as 
$\,  \frac{1}{j} \, \widehat{\varphi}_1 \,$. We smooth out 
$\, \widehat{\varphi}_j \,$ in a neighborhood of \f{0} with size of order
$\frac{1}{j}$ and obtain smooth functions 
$\, \varphi_j : \R \rightarrow \R \,$  with :
\begin{itemize}
\item	$\, 0 \le \varphi_j \le  \frac{1}{j} \,$ on \f{\R},
\item	$\, | \varphi^{\prime}_{j} | \le  \frac{2}{j} \,$ on \f{\R},
\item	$\, \varphi^{\prime\prime}_{j} \le  \frac{C_3}{j} \,$ on \f{\R},
\item	$\, \varphi_j \equiv 0 \,$ on $\, \R \setminus (- 1, 1) \,$,
\item	$\, \varphi^{\prime\prime}_{j} \, (0) \le - j \,$.
\end{itemize}

\begin{center}
\includegraphics[width=12cm]{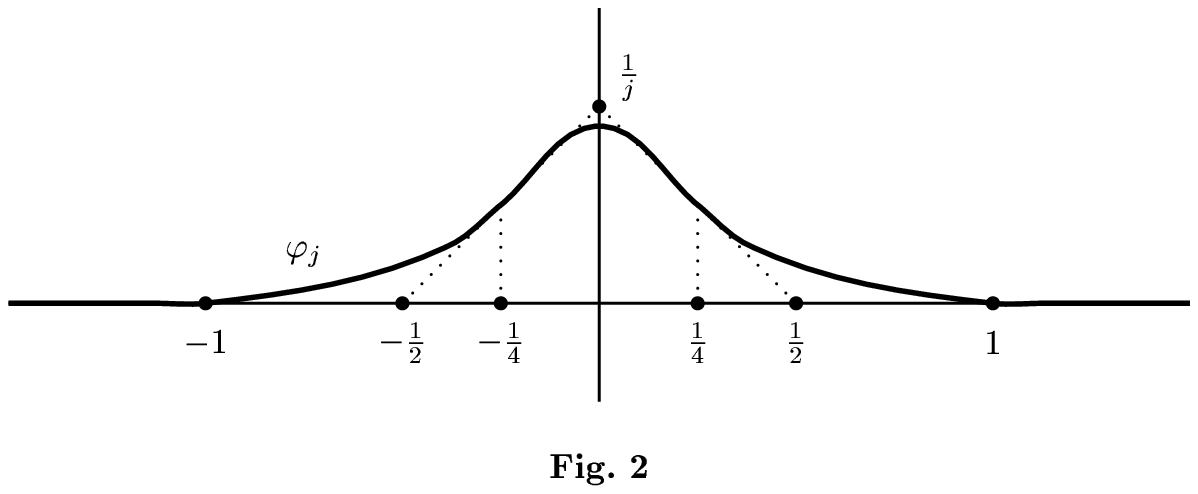}
\end{center}

One checks easily that $\, u_j := \varphi_j \circ \rho \,$ has all 
the required properties (using the formula 
$\, \Delta u_j = (\varphi^{\prime}_{j} \circ \rho)
\, \Delta \rho 
- (\varphi^{\prime\prime}_{j} \circ \rho) \, | d \rho |^2 \,$).

It remains to handle the case when the normal bundle of \f{N} in 
\f{M} has no trivial one-dimensional sub-bundle. Locally it always has. 
Hence we can cover \f{N} by open sets $\, V_k, k = 1, \ldots, \ell \,$, 
and find smooth functions $\, \rho_k : M \rightarrow \R \,$ having all 
the properties of the previous \f{\rho} on $\, V_k \,$ instead of \f{N}. 
Let $\, \chi_1, \ldots, \chi_{\ell + 1} \,$ be a partition of unity on 
\f{M} subordinate to the open covering 
$\, \exp \, (\pi^{- 1} (V_1) \cap \widetilde{U}), \ldots
\,,\; \exp \, (\pi^{- 1} (V_{\ell}) \cap \widetilde{U}) \,,\; 
M \setminus N \,$.  Then $\, \chi_1 |_{N}, \ldots, \chi_{\ell} |_{N} \,$ 
form a partition of unity on \f{N} subordinate to the covering 
$\, V_1, \ldots, V_{\ell} \,$. Again, it is easy to check that
$$
u_j \,:=\, \sum^{\ell + 1}_{k = 1} \, \chi_{k} \cdot
(\varphi_j \circ \rho_k)
$$
has the required properties.
\end{proof}

      \subsection{Distribution of $\mathbf{L^2}$-norm of eigenspinors and 
                  scalar curvature}

Next we show that eigenspinors try to avoid large scalar curvature,
that is, most of their $L^2$-norm is supported in the region of the 
manifold where the scalar curvature is small.

\begin{lemma}
\label{distributionlemma}
Let \f{M} be a closed Riemannian spin manifold.
Let $\, \Lambda, S_0, S_1 \in \R \,,\; S_0 < S_1 \,$. Suppose the
scalar curvature of \f{M} satisfies $\, {\rm scal} \ge S_0 \,$.
Put
$$
M_{+} := \{ x \in M\ |\ {\rm scal} \, (x) \ge S_1 \} .
$$
Then for any spinor satisfying
$$
\Vert D \psi \Vert^2_{L^2 (M)} \le \Lambda^2 \, \Vert \psi \Vert^2_{L^2 (M)}
$$
the following inequality holds:
$$
\int_{M_{+}} \, | \psi |^2 {\rm dvol} \le \frac{4 \Lambda^2 - S_0}{S_1 - S_0}
\, \int_{M} \, | \psi |^2 \, {\rm dvol} \;.
$$
\end{lemma}

\begin{proof}
Using the Lichnerowicz formula $\, D^2 = \nabla^{\ast}
\nabla +  \frac{1}{4}{\rm scal} \,$ we obtain
\begin{eqnarray*}
0 & \le &  \int_{M} \, | \nabla \, \psi |^2 \, {\rm dvol} \,=\,
 \int_{M} \langle \nabla^{\ast} \nabla \psi, \psi \rangle \, {\rm dvol} \\
& = &  \int_{M} \, \langle D^2 \psi, \psi \rangle \, {\rm dvol} 
- \frac{1}{4} \,  \int_{M} \, {\rm scal} \, | \psi |^2 \, {\rm dvol} \\
& \le & \Lambda^2 \,  \int_{M} \, | \psi |^2 \, {\rm dvol} \, 
-\,  \frac{S_0}{4} \,
 \int_{M \setminus M_{+}} \, | \psi |^2 \, {\rm dvol} \, -
 \frac{S_1}{4} \,  \int_{M_{+}} \, | \psi |^2 \, {\rm dvol} \;.
\end{eqnarray*}
Hence
$$
\left( \frac{S_1}{4} - \Lambda^2 \right) \, \int_{M_{+}} \, | \psi |^2 \, 
{\rm dvol}
\le \left( \Lambda^2 - \frac{S_0}{4} \right) \, 
\int_{M \setminus M_{+}} \, | \psi |^2 \, {\rm dvol}
$$
and therefore
$$
\left( \frac{S_1}{4} - \frac{S_0}{4} \right) \, 
\int_{M_{+}} \, | \psi |^2 \, {\rm dvol} \le
\left( \Lambda^2 - \frac{S_0}{4} \right) \, 
\int_{M} \, | \psi |^2 \, {\rm dvol} \;
$$
which is the claim of the lemma.
\end{proof}

\Remark{
The lemma can be slightly improved to
\begin{equation} \label{distribution-improved}
\int_{M_{+}} \, | \psi |^2 \, {\rm dvol} \le
\frac{4 \, \frac{n - 1}{n} \, \Lambda^2 - S_0}{S_1 - S_0}
\, \int_{M} \, | \psi |^2 \, {\rm dvol}
\end{equation}
with essentially the same proof. 
One combines the Lichnerowicz formula with the identity
$$
|\nabla\psi|^2 = |P \psi|^2 + \frac{1}{n}\,|D\psi|^2
$$ 
where $P$ is the {\it twistor operator}.
Estimating $\int_{M} \, | P \, \psi |^2 \, {\rm dvol}$ (instead
of $\int_{M} \, | \nabla \, \psi |^2 \, {\rm dvol}$) by $0$  one looses a 
little less.
Note that the estimate is nontrivial only for
$$
\Lambda^2 \le \frac{n}{n - 1} \, \frac{S_1}{4}
$$
because otherwise the coefficient on the right hand side of
(\ref{distribution-improved}) is greater than $1$. The estimate also 
shows that if $\, \Lambda^2 \le  \frac{n}{n - 1} \,  \frac{S_0}{4} \,$,
then $\, \psi = 0 \,$. This way we recover Friedrich's eigenvalue 
estimate \cite{friedrich80a}:
$$
\lambda^2 \ge \frac{n}{n - 1} \, \frac{S_0}{4} \;.
$$
In our application however $\, S_0 \,$ will typically be negative.
}

        \subsection{$\mathbf{L^2}$-norm on annular regions}

In the following we will look at the \geese{distance spheres} of a
compact submanifold \f{N} of a Riemannian manifold \f{M} given by
$$
S_N(r) \,:=\, \{ x \in M \, | \, {\rm dist} \, (x, N) \,=\, r \}
$$
for $\, r > 0 \,$. 
For $\, 0 \le R_1 < R_2 \,$ we define the annular region
$$
A_N(R_1, R_2) \,:=\, \{ x \in M\ |\ R_1 \le {\rm dist} (x, N) \le R_2 \}
\,=\, \bigcup_{R_1 \le r \le R_2} S_N(r) \;.
$$
Let \f{\nu} denote the unit normal vector field 
of $S_N(r)$ pointing away from \f{N}, let \f{W = -\nabla\nu} denote the 
Weingarten map (or shape operator) and let \f{H = 
\frac{1}{n-1}\,\mathrm{tr}(W)} be the mean curvature of $\, S_N(r) \,$ 
with respect to \f{\nu}. 
The following lemma tells us that under a suitable condition only
little $L^2$-norm of a spinor is contained in the annulus
$A_N(r,2r)$ compared to the bigger annulus $A_N(r,(2r)^{1/11})$
when $r$ is small.

\begin{lemma}
\label{l2annulus}
Let $M$ be an $n$-dimensional Riemannian spin manifold and let
$N \subset M$ be a compact submanifold of codimension $k \ge 3$.
Then there exists $0<R<1$ depending on $M$ and $N$ such that for any 
$0 < r \le \frac{1}{2}R^{11}$ and any smooth spinor $\varphi$ defined on 
$A_N(r,(2r)^{1/11})$ the following estimate holds
$$
\frac{\|\varphi\|^2_{L^2(A_N(r,2r))}}{\|\varphi\|^2_{L^2(A_N(r,(2r)^{1/11}))}}
\,\,\,\,\le\,\,\,\, 
10\,\, r^{5/2}
$$
provided
$$
\Re\int_{S_N(\rho)}\langle\nabla_\nu\varphi,\varphi\rangle\, dA \ge 0
$$
holds for all $\rho\in [r,(2r)^{1/11}]$.
\end{lemma}

\begin{proof}
From
$$
\frac{d}{d\rho} \, dA \,=\, - (n - 1) \, H \, dA
$$
and
$$
\frac{d}{d\rho} \, |\varphi|^2 \,=\, 2 \Re \, 
\langle \nabla_{\nu} \varphi, \varphi \rangle
$$
we conclude
\begin{eqnarray}
\frac{d}{d\rho} \int_{S_N(\rho)}|\varphi|^2\, dA
&=&
\Re \int_{S_N(\rho)} \{2\langle \nabla_{\nu} \varphi, \varphi \rangle 
- (n - 1) \, H \, |\varphi|^2\}\, dA \nonumber\\
&\ge&
- (n - 1) \, \int_{S_N(\rho)} H\,  |\varphi|^2\, dA .
\label{l2e1}
\end{eqnarray}
For $\, p \in N \,$ and \f{\nu} a unit normal vector to 
\f{N} at \f{p} let $\, P_{\nu} \,$ denote the orthogonal projection
$$
P_{\nu} \,:\, \nu^{\perp} \longrightarrow \nu^{\perp} \cap T_p N^{\perp} \;.
$$
In particular, we have the orthogonal decomposition
$$
T_p M \,=\, T_p N \oplus \R \, \nu \oplus {\rm im} (P_{\nu}) \;.
$$
Let $\, P_{\nu, \rho} \,$ denote the parallel translate of $\, P_{\nu} \,$
along $\, c (\rho) = \exp (\rho \, \nu) \,$. Then \f{W} is of the form
$$
W \, (c (\rho)) \,=\, - \frac{1}{\rho} \, P_{\nu, \rho} + C \, (\rho) \;,
$$
where $\, C \, (\rho) \,$ extends smoothly to $\, \rho = 0 \,$, compare
\cite{eschenburg-heintze90a}.
Since the rank of $\, P_{\nu} \,$ is equal to $\, k - 1 \,$ we conclude that
$$
(n-1)H = -\frac{k-1}{\rho} + \mbox{O}_{M, N} \, \left( 1 \right).
$$
For $0 < \rho \le R$, where $R$ is sufficiently small, we have
$$
\left|\mbox{O}_{M, N} \, \left( 1 \right)\right| \le \frac{1}{4\rho}
$$
and therefore
$$
(n-1)H \le -\frac{k-5/4}{\rho}.
$$
Plugging this into (\ref{l2e1}) we get 
$$
\frac{d}{d\rho} \int_{S_N(\rho)}|\varphi|^2\, dA
\ge
\frac{k-5/4}{\rho}\int_{S_N(\rho)}|\varphi|^2\, dA
\ge
\frac{7/4}{\rho}\int_{S_N(\rho)}|\varphi|^2\, dA.
$$
since $k\ge 3$.
Therefore
$$
\frac{d}{d\rho} \ln \int_{S_N(\rho)}|\varphi|^2\, dA \ge \frac{7/4}{\rho}
$$
and hence for $r_1 \le r_2$
\begin{eqnarray*}
\ln \left(\frac{\int_{S_N(r_2)}|\varphi|^2\, dA}
{\int_{S_N(r_1)}|\varphi|^2\, dA}\right)
&=&
\ln\left(\int_{S_N(r_2)}|\varphi|^2\, dA\right) -
\ln\left(\int_{S_N(r_1)}|\varphi|^2\, dA\right) \\
&\ge&
\frac{7}{4}\int_{r_1}^{r_2} \frac{d\rho}{\rho} \\
&=&
\frac{7}{4}\ln\left(\frac{r_2}{r_1}\right) .
\end{eqnarray*}
Exponentiation yields
\begin{equation}
\frac{\int_{S_N(r_2)}|\varphi|^2\, dA}
{\int_{S_N(r_1)}|\varphi|^2\, dA} 
\ge
\left(\frac{r_2}{r_1}\right)^{7/4}.
\label{l2e2}
\end{equation}
From (\ref{l2e1}) and $H < 0$ it follows that the function $\rho \mapsto 
\int_{S_N(\rho)}|\varphi|^2\, dA$ is monotonically increasing, so we have on
the one hand
\begin{equation}
\|\varphi\|^2_{L^2(A_N(r,2r))} = 
\int_r^{2r}\int_{S_N(\rho)}|\varphi|^2\, dA\, d\rho \le
r\, \int_{S_N(2r)}|\varphi|^2\, dA.
\label{l2e3}
\end{equation}
On the other hand by (\ref{l2e2}) with $r_1 = 2r$ and $r_2 = \rho$,
\begin{eqnarray}
\|\varphi\|^2_{L^2(A_N(2r,(2r)^{1/11}))} 
&=& 
\int_{2r}^{(2r)^{1/11}}\int_{S_N(\rho)}|\varphi|^2\, dA\, d\rho  \nonumber\\
&\ge&
\int_{S_N(2r)}|\varphi|^2\, dA 
\int_{2r}^{(2r)^{1/11}}\left(\frac{\rho}{2r}\right)^{7/4} \, d\rho\nonumber\\
&=&
\int_{S_N(2r)}|\varphi|^2\, dA \cdot
\frac{8r}{11}\left[\left(\frac{1}{2r}\right)^{5/2} - 1\right] .
\label{l2e4}
\end{eqnarray}
Combining (\ref{l2e3}) and (\ref{l2e4}) we obtain
\begin{eqnarray*}
\frac{\|\varphi\|^2_{L^2(A_N(r,(2r)^{1/11}))}}{\|\varphi\|^2_{L^2(A_N(r,2r))}}
&=&
1 + \frac{\|\varphi\|^2_{L^2(A_N(2r,(2r)^{1/11}))}}{\|\varphi\|^2_{L^2(A_N(r,2r))}} \\
&\ge&
1 + \frac{8}{11}\left[\left(\frac{1}{2r}\right)^{5/2} - 1\right] \\
&>&
\frac{1}{10}\left(\frac{1}{r}\right)^{5/2},
\end{eqnarray*}
and the lemma follows.
\end{proof}

	\subsection{Proof of Theorem~\ref{thm-spectralsurgery}}

By decreasing \f{\varepsilon} if necessary we may assume that
$\, (\Lambda - \varepsilon, \Lambda + \varepsilon) \,$ and
$\, (- \Lambda - \varepsilon, - \Lambda + \varepsilon) \,$ contain no
eigenvalues of $\, D_g \,$ and that any two distinct eigenvalues in
$\, (- \Lambda, \Lambda) \,$ are at distance at least $\, 2 \varepsilon \,$.
Write $\, \lambda_1 \le \lambda_2 \le \cdots \le \lambda_m \,$ for the
eigenvalues of $\, D_g \,$ in $\, (- \Lambda, \Lambda) \,$, each eigenvalue
being repeated according to its multiplicity. Denote the distance tube of 
radius \f{r} about \f{N} by $\, U_N(r) \,$, that is
$$
U_N(r) \,=\, \{ x \in M\ |\ {\rm dist} \, (x, N) < r \} \;.
$$
We will first show that there exists $\, R > 0 \,$ such that the Dirac 
operator $\, D_{\widetilde{g}} \,$ on any closed Riemannian spin manifold 
$\, (\widetilde{M},\widetilde{g}) \,$ containing an isometric copy of 
$\, M \setminus U_N(r) \,$ (together with its spin structure) has 
eigenvalues
$\mu_1 \le \mu_2 \le \cdots \le \mu_{m}$
with $\, | \lambda_i - \mu_i | < \varepsilon \,$ provided $\, r < R \,$.
It may however have further eigenvalues in the interval 
$\, (- \Lambda, \Lambda) \,$. This will be excluded later under 
additional assumptions.

To prove the existence of the eigenvalues $\, \mu_{i}$ let
$\chi_{r} : M \rightarrow \R$
be a smooth cut--off function with the following properties
\begin{itemize}
\item	$\, 0 \le \chi_r \le 1 \,$ on \f{M},
\item	$\, \chi_r \equiv 0 \,$ on $\, U_N(r) \,$,
\item	$\, \chi_r \equiv 1 \,$ on $\, M \setminus U_N(2 r) \,$,
\item	$\, | \nabla \chi_r | \le  \frac{2}{r} \,$ on \f{M}.
\end{itemize}
Since all norms on a finite dimensional vector space are equivalent
there exists a constant $\, C_1 > 0 \,$ such that
$$
\Vert \varphi \Vert_{L^{\infty} (M)} \le C_1 \,
\Vert \varphi \Vert_{L^2 (M)}
$$
for all eigenspinors \f{\varphi} of $\, D_g \,$ with eigenvalue in
$\, (- \Lambda, \Lambda) \,$. Moreover, there is a constant 
$\, C_2 > 0 \,$ such that
$$
{\rm vol} \, (U_N(2 r)) \le C_2 \, r^k \;.
$$
Let \f{\varphi_i} be an eigenspinor on \f{M} for the eigenvalue 
$\, \lambda_i \,$. Since $\, \chi_r \, \varphi_i \,$ has its support 
in $\, M \setminus U_N(r) \,$ it can also be regarded as a spinor on 
$\, \widetilde{M} \,$. We plug it into the Rayleigh quotient of 
$\, D_{\widetilde{g}} - \lambda_i \,$ :
\begin{eqnarray*}
 \frac{ \Vert (D_{\widetilde{g}} - \lambda_i) \,
(\chi_r \, \varphi_i) \Vert^{2}_{L^2 (\widetilde{M})}}
{\Vert \chi_r \, \varphi_i \Vert_{L^2 (\widetilde{M})}^2} &=&
 \frac{\Vert \nabla \chi_r \cdot \varphi_i \Vert^2_{L^2 (U_N(2r))}}
{\Vert \chi_r \, \varphi_i \Vert_{L^2 (M)}^2} \\
&\le&  \frac{ \frac{4}{r^2} \, \Vert \varphi_i \Vert^2_{L^2 (U_N(2r))}}
{\Vert \varphi_i \Vert^2_{L^2 (M)} 
- \Vert \varphi_i \Vert^2_{L^2 (U_N(2r))}} \;.
\end{eqnarray*}
From
\begin{eqnarray*}
\Vert \varphi_i \Vert^2_{L^2 (U_N(2r))} & \le &
\Vert \varphi_i \Vert^2_{L^{\infty}(U_N(2r))} \, {\rm vol} \, (U_N(2r)) \\
& \le & C^2_1 \, \Vert \varphi_i \Vert^2_{L^2 (M)} \, C_2 \, r^k
\end{eqnarray*}
we conclude
\begin{eqnarray*}
 \frac{\Vert (D_{\widetilde{g}} - \lambda_i) \, (\chi_r \,
\varphi_i) \Vert^2_{L^2 (\widetilde{M})}}
{\Vert \chi_r \, \varphi_i \Vert^2_{L^2 (\widetilde{M})}}
& \le &
\frac{ \frac{4}{r^2} \, C^2_1 \, C_2 \, r^k }
{1 - C^2_1 \, C_2 \, r^k}
\\
& = &  \frac{4 \, C^2_1 \, C_2 \, r^{k - 2}}
{1 - C_1^2 \, C_2 \, r^k} \\
& < & \varepsilon^2
\end{eqnarray*}
for \f{r} sufficiently small because $\, k \ge 3 \,$. Therefore the 
operator $\, D_{\widetilde{g}} - \lambda_i \,$ has an eigenvalue in 
the interval $\, (- \varepsilon, \varepsilon) \,$, hence 
$\, D_{\widetilde{g}} \,$ has an eigenvalue $\, \mu_i \,$ in 
$\, (\lambda_i - \varepsilon, \lambda_i + \varepsilon) \,$.
If $\, \lambda_i \,$ has multiplicity \f{\ell} this yields 
\f{\ell} eigenvalues of $\, D_{\widetilde{g}} \,$ in 
$\, (\lambda_i - \varepsilon, \lambda_i + \varepsilon) \,$.

It remains to be shown that in the situation of the theorem 
$\widetilde{g}$ can be chosen so that the operator
$\, D_{\widetilde{g}} \,$ has no more than \f{m} eigenvalues in
$\, (- \Lambda, \Lambda) \,$. According to 
Proposition~\ref{prop-incrscalcurv} there is a metric on \f{M} 
arbitrarily close to \f{g} in the $C^1$-topology such that with 
respect to this metric
\begin{itemize}
\item	$\, {\rm scal} \ge S_0 \,$ on all of $M$,
\item	$\, {\rm scal} \ge 2 S_1 \,$ on a neighborhood
	$\, U_0 \,$ of \f{N} with $S_1$ arbitrarily large.
\end{itemize}

Since the eigenvalues of $\, D_g \,$ depend continuously on the 
Riemannian metric with respect to the $C^1$-topology 
(see for example \cite[Prop.~7.1]{baer96b}) we may assume without loss 
of generality that the scalar curvature of \f{g} itself has these 
properties with $\, S_1 \,$ so large that
$$
\frac{S_1 - 4 \Lambda^2}{S_1 - S_0} \ge \left(
\frac{\Lambda + \varepsilon/2}{\Lambda + \varepsilon} \right)^2 \;.
$$
Now choose $\, r > 0 \,$ so small that
\begin{itemize}
\item	$\, r <  \frac{\varepsilon^4}{2^8\cdot 100 \cdot (m + 1)^4} \,$,
\item	$\, U_N((2r)^{1/11}) \subset U_0 \,$,
\item	$\, (2r)^{1/11} \,$ is no larger than the \f{R} in 
Lemma~\ref{l2annulus}.
\end{itemize}

\begin{center}
\includegraphics[width=12cm]{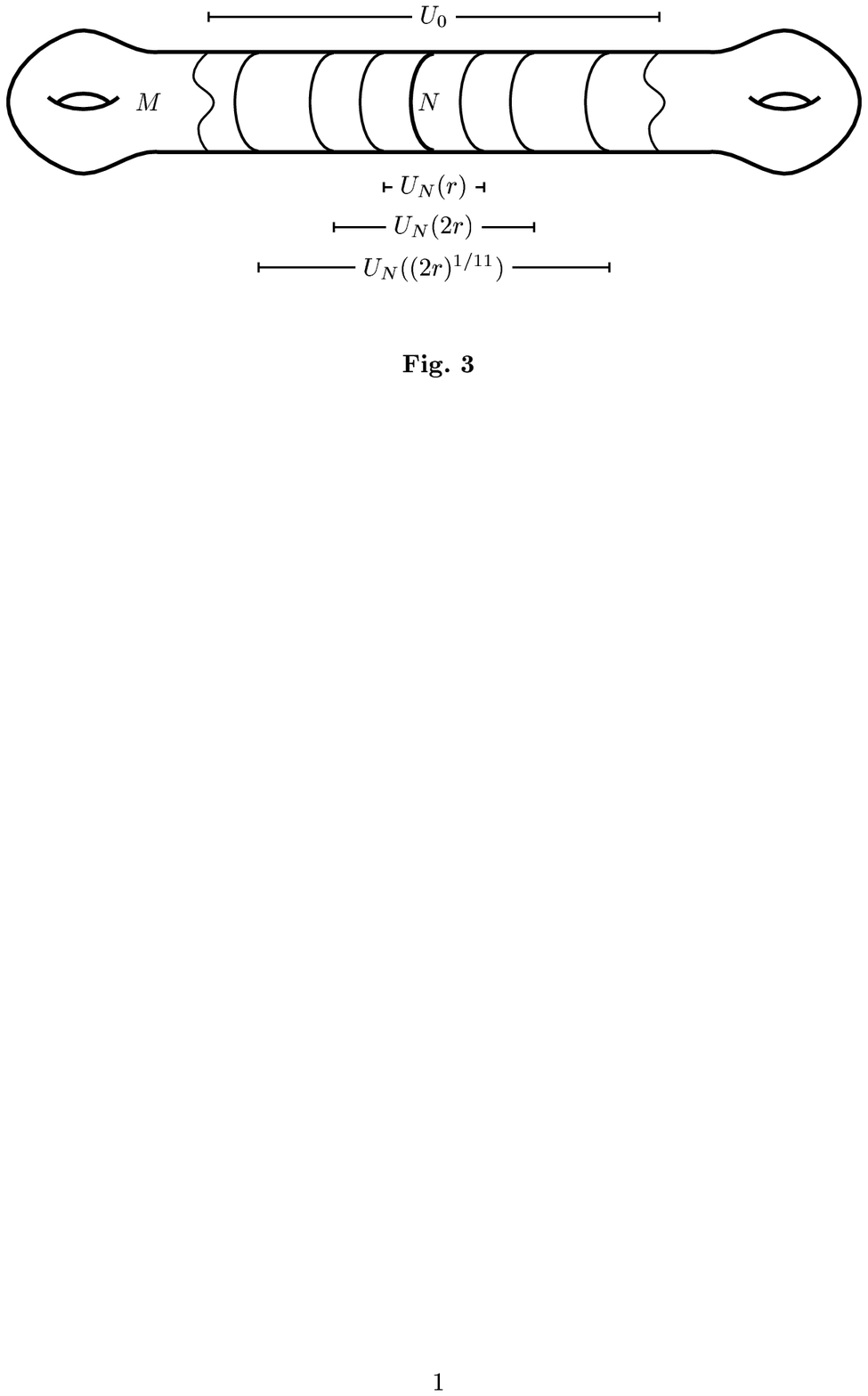}
\end{center}

We perform the surgery along $N$ in the neighborhood $\, U = U_N(r) \,$.
Hence $\widetilde{M}$ is of the form $\widetilde{M} = (M \setminus U_N(r))
\cup \widetilde{U}$.
Surgery in codimension $\, \ge 3 \,$ does not decrease scalar curvature 
too much if the metric $\, \widetilde{g} \,$ on $\, \widetilde{M} \,$ 
is chosen properly, see \cite[Proof of Theorem~A]{gromov-lawson80a} and
\cite[Proof of Theorem~3.1]{rosenberg-stolz98a}. We may assume
\begin{itemize}
\item	$\, {\rm scal}_{\widetilde{g}} \ge S_0 \,$ on all of 
	$\, \widetilde{M} \,$,
\item	$\, {\rm scal}_{\widetilde{g}} \ge S_1 \,$ on 
	$\, \widetilde{U} \,$.
\end{itemize}

We will show that with these choices of $\, r \,$ and 
$\, \widetilde{g} \,$ the Dirac operator of $\, \widetilde{M} \,$ can 
have no more than \f{m} eigenvalues in $\, (- \Lambda, \Lambda) \,$. 
Assume the contrary, that is, assume there is an $(m + 1)$-dimensional 
space $\, \widetilde{{\mathcal H}} \,$ of spinors on 
$\, \widetilde{M} \,$ spanned by eigenspinors with eigenvalues in 
$\, (- \Lambda, \Lambda) \,$. The function $\, \chi_{r} \,$ can also 
be considered as a cut--off function on $\, \widetilde{M} \,$ since its 
support is contained in $\, M \setminus U_N(r) \,$ which is 
contained in $\, \widetilde{M} \,$.  
For every spinor \f{\varphi} on $\, \widetilde{M} \,$ the spinor 
$\, \chi_{r} \, \varphi \,$ can be considered a spinor also on \f{M}. 
The space
$$
{\mathcal H} \,:=\, \{ \chi_{r} \, \varphi \, | \,
\varphi \in \widetilde{{\mathcal H}} \}
$$
has the same dimension $\, m + 1 \,$ as $\widetilde{{\mathcal H}}$ by the 
unique-continuation property of the Dirac operator (see for example 
\cite{booss-wojciechowski93a}) 
and we consider it as a space of spinors on \f{M}. We will show that 
$$
\frac{\Vert D_g \psi \Vert^2_{L^2 (M)}}{\Vert \psi \Vert^2_{L^2 (M)}}
< \left( \Lambda + \varepsilon \right)^2 \;
$$
for all nonzero $\, \psi \in {\mathcal H}$.
Then $\, D_g \,$ must have at least $\, m + 1 \,$ eigenvalues in
$\, (- \Lambda - \varepsilon, \Lambda + \varepsilon) \,$
which is a contradiction.

Let $\, \psi := \chi_{r} \, \varphi \in {\mathcal H} \,$. 
We write $\, \varphi = \varphi_1 + \cdots + \varphi_{m + 1} \,$ with 
$D_{\widetilde{g}}$-eigenspinors $\, \varphi_j \,$ (some possibly zero)
for the eigenvalues $\mu_j\in(-\Lambda,\Lambda)$.
By the assumption on the scalar curvature of $\, \widetilde{g} \,$ 
(recall that $\scal \ge S_1$ on $\widetilde{U}$ and $\scal \ge 2\, S_1$ on $A_N(r,2r)$)
and by Lemma~\ref{distributionlemma} applied to $\, \widetilde{M} \,$ we 
have
$$
\Vert \varphi \Vert^2_{L^2 (\widetilde{U} \cup A_N(r,2r))} 
\le
\frac{4 \Lambda^2 - S_0}{S_1 - S_0} \, 
\Vert \varphi \Vert^2_{L^2 (\widetilde{M})} \;.
$$
Thus
\begin{eqnarray*}
\Vert \psi \Vert^2_{L^2 (M)} 
& = & 
\Vert \chi_{r} \, \varphi \Vert^2_{L^2 (M)} 
\ge
\Vert \varphi \Vert^2_{L^2 (M \setminus U_N(2 r))} \\
& \ge & 
\left[ 1 -  \frac{4 \Lambda^2 - S_0}{S_1 - S_0} \right]
\cdot \Vert \varphi \Vert^2_{L^2 (\widetilde{M})} \\
& = &  
\frac{S_1 - 4 \Lambda^2}{S_1 - S_0} \, \Vert \varphi
\Vert^2_{L^2 (\widetilde{M})} \\
& \ge & 
\left(  \frac{\Lambda +  
\frac{\varepsilon}{2}}{\Lambda + \varepsilon} \right)^2 \, 
\Vert \varphi \Vert^2_{L^2 (\widetilde{M})} \;.
\end{eqnarray*}
Secondly,
\begin{eqnarray}
\Vert D_g \, \psi \Vert_{L^2 (M)} & = & \Vert D_{\widetilde{g}}
\, (\chi_{r} \, \varphi) \Vert_{L^2 (\widetilde{M})} \nonumber\\
& = & \Vert \nabla \chi_{r} \cdot \varphi + \chi_{r} \, D_{\widetilde{g}} \,
\varphi \Vert_{L^2 (\widetilde{M})} \nonumber\\
& \le & \Vert \nabla \chi_{r} \cdot \varphi \Vert_{L^2 (\widetilde{M})}
+ \Vert D_{\widetilde{g}} \, \varphi \Vert_{L^2 (\widetilde{M})} \nonumber\\
& \le &  \frac{2}{r} \, \Vert \varphi \Vert_{L^2 (A_N(r, 2r))}
+ \Lambda \, \Vert \varphi \Vert_{L^2 (\widetilde{M})} \;.
\label{eganzneu}
\end{eqnarray}
Now we justify that we can apply Lemma~\ref{l2annulus} to the eigenspinors
$\varphi_j$.
Fix $\rho\in [r,(2r)^{1/11}]$.
We have to show
$$
\Re\int_{S_N(\rho)}\langle\nabla_\nu\varphi_j,\varphi_j\rangle\, dA \ge 0.
$$
Set $\Mhat := \widetilde{U} \cup A_N(r,\rho)$.
Then $\Mhat \subset \widetilde{M}$ is a compact manifold with boundary
$\partial\Mhat = S_N(\rho)$ and $\scal \ge S_1$ on $\Mhat$.

\begin{center}
\includegraphics[width=12cm]{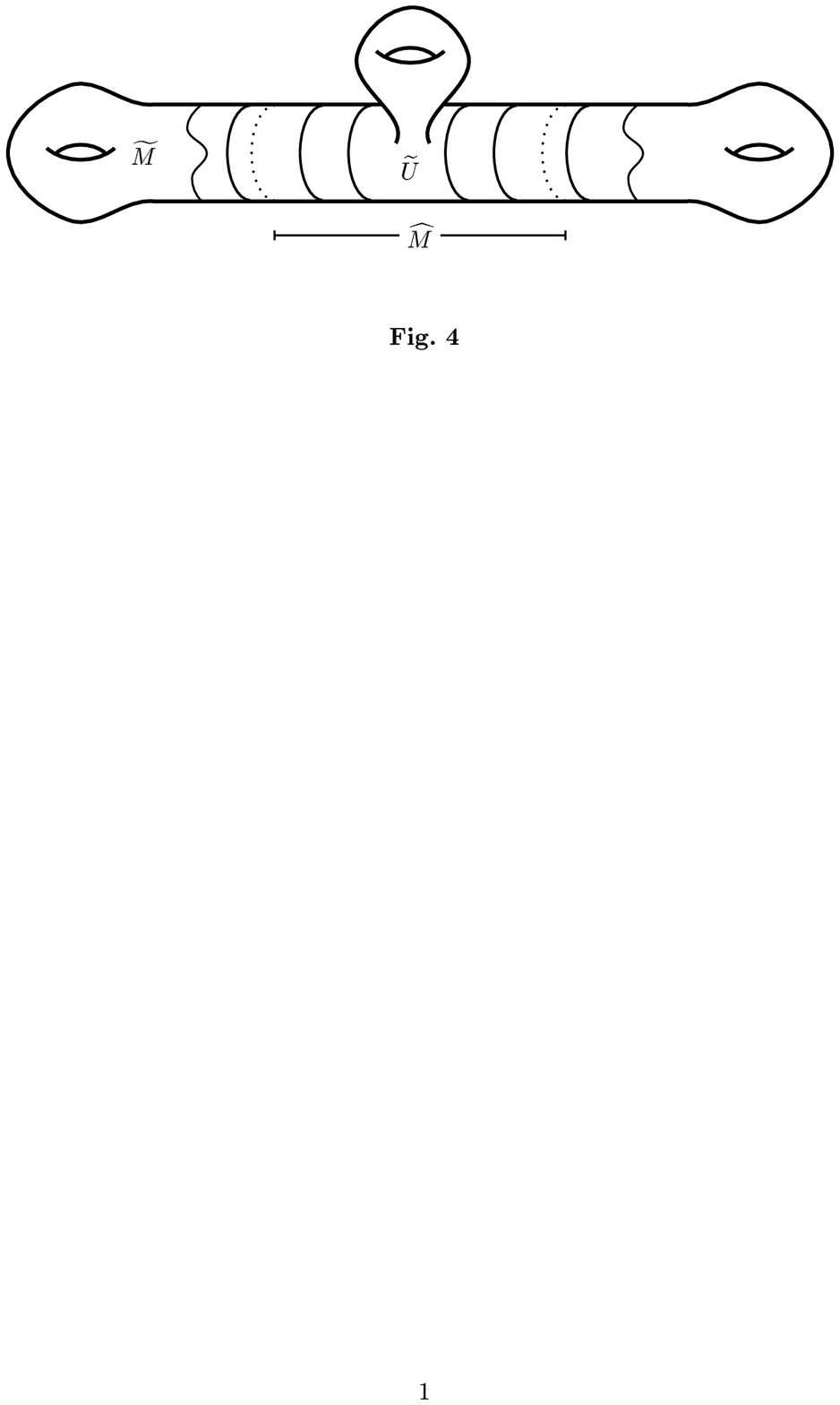}
\end{center}

The Lichnerowicz formula and a partial integration yield
\begin{eqnarray*}
\mu_j^2 \, \|\varphi_j\|_{L^2(\Mhat)}^2 
&=&
\int_{\Mhat} \langle D^2\varphi_j,\varphi_j\rangle\, \dvol \\
&=&
\int_{\Mhat} \langle \nabla^\ast\nabla\varphi_j,\varphi_j\rangle\, \dvol
+ \frac{1}{4}\int_{\Mhat} \scal\, |\varphi_j|^2 \, \dvol \\
&\ge&
\|\nabla \varphi_j\|_{L^2(\Mhat)}^2 
- \int_{\partial\Mhat} \langle\nabla_\nu\varphi_j,\varphi_j\rangle\, dA
+ \frac{S_1}{4}\, \|\varphi_j\|_{L^2(\Mhat)}^2 .
\end{eqnarray*}
Hence 
$$
\int_{S_N(\rho)}\langle\nabla_\nu\varphi_j,\varphi_j\rangle\, dA \ge 
\left(\frac{S_1}{4}-\Lambda^2\right)\, \|\varphi_j\|_{L^2(\Mhat)}^2 \ge 0.
$$
Thus Lemma~\ref{l2annulus} can be applied and yields
\begin{eqnarray*}
\Vert \varphi_j \Vert^2_{L^2 (A_N(r, 2r))}
& \le &  
10 \, r^{5/2} \, \Vert \varphi_j \Vert^2_{L^2 (A_N(r, (2r)^{1/11}))} \\
& \le & 
10 \, r^{5/2} \, \Vert \varphi_j \Vert^2_{L^2 (\widetilde{M})} \\
& \le & 
10 \, r^{5/2} \, \Vert \varphi \Vert^2_{L^2 (\widetilde{M})} \;.
\end{eqnarray*}
Thus
\begin{eqnarray*}
\frac{2}{r} \, \Vert \varphi \Vert_{L^2 (A_N(r, 2r))}
& \le &  \frac{2}{r} \, \left\{ \Vert \varphi_1 
\Vert_{L^2 (A_N(r, 2r))}
+ \cdots + \Vert \varphi_{m + 1} \Vert_{L^2 (A_N(r, 2r))} \right\} \\
& \le &  \frac{2}{r} \, (m + 1) \, \sqrt{10 \cdot r^{5/2}} \,
\Vert \varphi \Vert_{L^2 (\widetilde{M})} \\
& = & 2 \, \sqrt{10} \, (m + 1) \, r^{1/4} \,
\Vert \varphi \Vert_{L^2 (\widetilde{M})} \\
& < & \frac{\varepsilon}{2} \, 
\Vert \varphi \Vert_{L^2 (\widetilde{M})} \;.
\end{eqnarray*}
Plugging this into (\ref{eganzneu}) gives
$$
\Vert D_g \psi \Vert_{L^2 (M)} < 
\left( \Lambda + \frac{\varepsilon}{2} \right)
\, \Vert \varphi \Vert_{L^2 (\widetilde{M})}
$$
and therefore
$$
\frac{\Vert D_g \psi \Vert^2_{L^2 (M)}}{\Vert \psi \Vert^2_{L^2 (M)}} <
\frac{ \left( \Lambda + \frac{\varepsilon}{2} \right)^2 \,
\Vert \varphi \Vert^2_{L^2 (\widetilde{M})}}{\left( \frac{\Lambda +
\frac{\varepsilon}{2}}{\Lambda + \varepsilon}
\right)^2 \, \Vert \varphi \Vert^2_{L^2 (\widetilde{M})} }
\,=\, (\Lambda + \varepsilon)^2 \;
$$
which completes the proof of Theorem~\ref{thm-spectralsurgery}.

\section{Generic metrics and harmonic spinors}
\label{sec-harmonic}

The {\em alpha-genus} is a ring homomorphism 
$\alpha_*: \Omega_*^{\text{spin}} \rightarrow KO^{-*}(\text{pt})$ 
from the spin bordism ring to the real K-theory ring. 
This means that the alpha-genus $\alpha(M)$ of a spin manifold $M$
depends only on the spin bordism class of $M$ and that it is additive with 
respect to disjoint union and multiplicative with respect to product 
of spin manifolds. Choosing generators we have 
$$
KO^{-n} (\text{pt}) = \left\{
\begin{array}{c l}
\Za 	& \quad \text{if } n \equiv 0,4 \mod 8, \\
\Za/2 	& \quad \text{if } n \equiv 1,2 \mod 8, \\
0 	& \quad \text{otherwise}.
\end{array}
\right.
$$
In this identification $\alpha_{8k} = \Ahat_{8k}$ and 
$\alpha_{8k+4} = \half \Ahat_{8k+4}$ where 
$\Ahat_*: \Omega_*^{\text{spin}} \rightarrow \Za$ 
is the $\Ahat$-genus. 

There is a generalization $\mathcal{D}$ of the Dirac operator called the 
{\em $Cl_n$-linear Atiyah-Singer operator} \cite{lawson-michelsohn89a}.
This operator acts on sections of a bundle associated to the spin 
structure with fiber the Clifford algebra $Cl_n$ and commutes with 
the right action of $Cl_n$ on this bundle. The operator $\mathcal{D}$ 
has a ``Clifford index'' $\ind_n(\mathcal{D}) \in KO^{-n}(\text{pt})$. 
This bundle of Clifford algebras decomposes into a sum of copies of the 
spinor bundle and the Clifford index is related to the index of the 
Dirac operator acting on sections of the spinor bundle as follows.
$$
\ind_n(\mathcal{D}) =  
\left\{
\begin{array}{c l}
\ind(D^+)		& \quad \text{if } n \equiv 0 \mod 8, \\
\half \ind(D^+)		& \quad	\text{if } n \equiv 4 \mod 8, \\
\dim \ker D \mod 2 	& \quad \text{if } n \equiv 1 \mod 8, \\
\dim \ker D^+ \mod 2	& \quad \text{if } n \equiv 2 \mod 8.
\end{array}
\right.
$$ 
The index theorem for $\mathcal{D}$ states that $\ind_n(\mathcal{D})$ 
coincides with $\alpha_n(M)$ and this gives a lower bound on the 
dimension of the kernel of the Dirac operator $D$ on a manifold $M$,
$$
\dim \ker D \geq  
\left\{
\begin{array}{c l}
|\Ahat(M)|	& \quad \text{if } n \equiv 0,4 \mod 8, \\
|\alpha(M)| 	& \quad \text{if } n \equiv 1 \mod 8, \\
2|\alpha(M)|  	& \quad \text{if } n \equiv 2 \mod 8.
\end{array}
\right.
$$ 
The left hand side of this inequality depends on the Riemannian 
metric while the right hand side is a differential topological invariant. 
In the following definition we give a name to those Riemannian metrics 
for which equality holds.

\Definition{
A Riemannian metric $g$ on a compact spin manifold $M$ is called 
{\em $D$-minimal} if the kernel of the Dirac operator $D_g$ is no larger 
than it is forced to be by the index theorem.
}

In dimensions $3,5,6,7 \mod 8$ the index of $D$ always vanishes and a metric
$g$ is $D$-minimal if $\ker D_g = 0$. 

Denote by $h_g$ the dimension of the kernel of the Dirac operator 
$D_g$ and if the dimension is even denote by $h_g^+, h_g^-$ the 
dimensions of the kernels of $D_g^+, D_g^-$. On an $n$-dimensional 
manifold $M$ with a $D$-minimal metric $g$ the numbers 
$h_g, h_g^+, h_g^-$ are determined by $\alpha$ according to the 
following table.

\bigskip
\begin{center}
\begin{tabular}{| c | c | c |}
\hline
$n \mod 8$ 	& $\alpha(M)$	& 					\\
\hline\hline
 & $\geq 0$	& $h_g^+ =\Ahat$, \quad $h_g^- = 0$ 			\\
\cline{2-3}
\raisebox{1.8ex}[-1.8ex]
{$0,4$}		& $ < 0$	& $h_g^+ = 0$,	\quad $h_g^- = -\Ahat$ 	\\
\hline
 		& $0$		& $h_g = 0$				\\
\cline{2-3}
\raisebox{1.8ex}[-1.8ex]
{$1$}		& $1$		& $h_g = 1$				\\
\hline
		& $0$		& $h_g^+ = h_g^- = 0$			\\
\cline{2-3}
\raisebox{1.8ex}[-1.8ex]
{$2$}		& $1$		& $h_g^+ = h_g^- = 1$			\\
\hline
$3,5,6,7$ 	& $0$		& $h_g = 0$				\\
\hline
\end{tabular}
\end{center}
\bigskip

With this terminology Conjecture~B from the introduction can be 
phrased as

{\bf Conjecture B.}
{\em 
On any compact connected spin manifold a generic metric is $D$-minimal. 
}

The first result in this direction is by Anghel \cite{anghel96a} who 
studies twisted Dirac operators on a manifold with a fixed Riemannian 
metric and shows that in dimensions $\leq 4$ a generic choice of 
connection on the twisting bundle gives a minimal kernel. The 
idea of the proof is to show that a non-minimal connection can be 
deformed to give nearby minimal connections. 
Such deformations are 
easily found only under the strong restriction on dimension. Using the 
formulas by Bourguignon and Gauduchon \cite{bourguignon-gauduchon92a} 
for the variation of the spinor bundle and the Dirac operator with the 
Riemannian metric Maier \cite{maier97a} adapts the method of Anghel to 
prove Conjecture~B in the same dimensions.

Denote by $\mathcal{R}(M)$ the space of smooth Riemannian metrics on a 
compact spin manifold $M$ and by $\mathcal{R}_{\rm{min}}(M)$ the subset 
of $D$-minimal metrics. We have not yet made clear what we mean by the 
term ``generic''. Standard results from perturbation theory give the 
following proposition (see \cite[Section {VII.1.3}]{kato76a},
\cite{anghel96a}, \cite{maier97a}).

\begin{prop} 
\label{opendense}
$\mathcal{R}_{\rm{min}}(M)$  is open in the $C^1$-topology on 
$\mathcal{R}(M)$ and if it is non-empty it is dense in all 
$C^k$-topologies, $k \geq 1$.
\end{prop}

We take ``generic'' to mean ``in a set which is open in the 
$C^1$-topology on $\mathcal{R}(M)$ and dense in all $C^k$-topologies, 
$k \geq 1$''. Denote by $\mmin$ the class of compact spin manifolds 
possessing at least one $D$-minimal metric. It follows from 
Proposition~\ref{opendense}
that Conjecture B is true precisely for manifolds in $\mmin$.

The rest of this paper is devoted to a study of the class $\mmin$.

	\subsection{Properties of the class $\mmin$}

Note that the conclusion of Conjecture B may not hold
if the manifold is not connected. Let $M$ be a manifold with 
$\alpha(M) \neq 0$ and let $-M$ be the same manifold with reversed 
orientation. Then the disjoint union $M + (-M)$ has 
$\alpha(M + (-M)) = 0$ but the kernel of the Dirac operator cannot be
trivial because it then would be trivial also on $M$. The following 
proposition tells us that we can take disjoint unions (and later 
connected sums) in the class $\mmin$ whenever there is no cancellation 
of the alpha-genera.

\begin{prop} \label{prop-disjointunion} 
Let $M_1,M_2 \in \mmin$ with $\dim M_1 = \dim M_2 = n$. Suppose
\begin{enumerate}
\item $n \equiv 0,4 \mod 8$ and $\alpha(M_1) \alpha(M_2) \geq 0$, or
\item $n \equiv 1,2 \mod 8$ and $\alpha(M_1) \alpha(M_2) =0$, or
\item $n \equiv 3,5,6,7 \mod 8$.
\end{enumerate}
Then the disjoint union $M_1 + M_2 \in \mmin$. 
\end{prop}

\begin{proof}
We only treat the first case $n = 4k$ and $\alpha(M_1) \alpha(M_2) \geq 0$.
The other cases are similar.

Let $g_1$ and $g_2$ be $D$-minimal metrics on $M_1$ and $M_2$ and let 
$g$ be the metric on $M_1 + M_2 $ given by $g_1$ and $g_2$. Suppose 
$\alpha(M_1 + M_2) = \alpha(M_1) + \alpha(M_2) =0$, then 
$\alpha(M_1) = \alpha(M_2) = 0$ since we assume 
$\alpha(M_1) \alpha(M_2) \geq 0$. This means that  $h_{g_1} = h_{g_2} = 0$ 
so $h_{g} = 0$ and $g$ is $D$-minimal. Next suppose 
$\alpha(M_1 + M_2) = \alpha(M_1) + \alpha(M_2)  > 0$. Then we cannot have  
$\alpha(M_1) \leq 0$ and $\alpha(M_2) \leq 0$. By assumption $\alpha(M_1)$ 
and $\alpha(M_2)$ do not have different signs and we must have 
$\alpha(M_1) \geq 0$ and $\alpha(M_2) \geq 0$. Since $g_1$ and $g_2$ are 
$D$-minimal we have $h_g^- =   h_{g_1}^- + h_{g_2}^- = 0 + 0 = 0$ so 
$g$ is $D$-minmimal. The case $\alpha(M_1 + M_2) < 0$ is analogous to 
$\alpha(M_1 + M_2) > 0$.   
\end{proof}

A metric $g$ with positive scalar curvature is $D$-minimal since the 
Lichnerowicz formula forces $h_g = 0$. This gives a rich source of 
manifolds in the class $\mmin$.

\begin{prop} \label{prop-lichnerowicz}
Let $M$ be a compact spin manifold which has a metric
of positive scalar curvature. Then $M \in \mmin$. 
\end{prop}

A manifold with a positive scalar curvature metric has vanishing 
alpha-genus \cite{hitchin74a}. We will need a set of manifolds in 
$\mmin$ whose alpha-genera generate $KO^{-*}(\text{pt})$.

\begin{prop}  \label{prop-generators}
For each $n \equiv 0,1,2,4 \mod 8$, $n \geq 1$, there is a $V_n \in \mmin$ 
with $\dim(V_n) = n$ and $\alpha(V_n) = 1$. 
\end{prop}

\begin{proof}
We prove this by induction starting with $n=1,2,4,8$. 

Let $V_1$ be the circle with the spin structure which is not spin 
bordant to zero, let $V_2 = V_1 \times V_1$, and let $V_4$ be the 
K3 surface. Let $V_8$ be a compact Riemannian $8$-manifold with 
holonomy $\text{Spin}(7)$ \cite{joyce96a}. For $V_1$ we have 
$\alpha(V_1)=1$ and $h =1$ so $V_1 \in \mmin$. By conformal invariance 
any metric $g$ on $V_2$ has $h_g^+ = h_g^- =1$ and since $\alpha(V_2) = 1$ 
we have $V_2 \in \mmin$. The manifold $V_4$ has 
$\Ahat(V_4) = 2\alpha(V_4) = 2$ so for any metric $g$ we have 
$h_g^+ \geq 2$. Let $g_0$ be a Ricci flat metric on $V_4$.
By the Lichnerowicz formula all harmonic spinors are parallel for this metric,
so $h_{g_0}^+ \leq \text{rank} (\Sigma^+ V_4) = 2$. This means that 
$h_{g_0}^+ = 2$ so $g_0$ is $D$-minimal and $V_4 \in \mmin$. The 
manifold $V_8$ with a $\text{Spin}(7)$-holonomy metric $g_8$ has a 
one-dimensional space of parallel spinors \cite{wang89a} and since 
it is Ricci-flat there are no other harmonic spinors. This means that 
$\alpha(V_8) = \Ahat(V_8) = 1 = h_{g_{_8}}^+$ so $g_8$ is $D$-minimal 
and $V_8 \in \mmin$. 

Next suppose that for some $n= 8k+1, 8k+2, 4k+4$, $k \geq 0$ we have 
a manifold $V_n$ with the desired properties.
Set $V_{n+8} := V_8 \times V_n$ with the product metric $g_{n+8} =
g_8 + g_n$, where $g_n$ is a $D$-minimal metric on $V_n$. Then 
$\alpha(V_{n+8}) = \alpha(V_8) \alpha (V_n) = 1$. The spinor bundle 
on $V_{n+8}$ is isomorphic to the outer tensor product of the spinor 
bundles on $V_8$ and $V_n$,   
$\Sigma V_{n+8} = \Sigma V_8 \widehat{\otimes} \Sigma V_n$, and the 
Dirac operators are related by 
$D_{g_{n+8}}^2 = D_{g_{_8}}^2 \otimes 1 + 1 \otimes D_{g_n}^2$. 
From this we see that the harmonic spinors on the product are given 
precisely by 
products of harmonic spinors on the factors. In particular we have 
$h_{g_{n+8}} = h_{g_{_8}} \cdot h_{g_n} = h_{g_n}$. Since $g_n$ is 
$D$-minimal $g_{n+8}$ is a $D$-minimal metric on $V_{n+8}$ and 
$V_{n+8} \in \mmin$. 
\end{proof}

The most important step in our study of the class $\mmin$ is the 
following application of the surgery theorem for the Dirac spectrum. 

\begin{prop} \label{prop-surgery}
$\mmin$ is closed under surgery in codimension $\geq 3$.
\end{prop}

\begin{proof}
Let $M \in \mmin$ and suppose $\tM$ is obtained from $M$ by surgery in 
codimension $\geq 3$. Let $g$ be a $D$-minimal metric on $M$ and choose 
$\Lambda > 0$ so that the interval $[-\Lambda, \Lambda]$ does not contain 
any non-zero eigenvalues of $D_g$. By Theorem \ref{thm-spectralsurgery} 
there is a metric $\tg$ on $\tM$ with the property that $D_{\tg}$ has 
precisely $h_g$ many eigenvalues in the interval $[-\Lambda, \Lambda]$ 
(counting multiplicity). Hence $h_{\tg} \le h_g$. Since surgery does not 
change the spin bordism class of a manifold we have 
$\alpha(\tM) = \alpha(M)$ and therefore $h_{\tg} \ge h_g$ because $g$ is 
$D$-minimal. We conclude that $h_{\tg} = h_g$ and therefore the metric $\tg$ 
is also $D$-minimal.
\end{proof}

\subsection{Simply connected manifolds}
We want to use Proposition \ref{prop-surgery} to construct new manifolds 
in the class $\mmin$. Without restrictions on the codimensions we can 
start with a spin manifold $N$ and do a sequence of surgeries ending up 
with a spin manifold $M$ if and only if $N$ and $M$ are spin bordant. 
The problem in our situation is that Proposition~\ref{prop-surgery} allows
only surgeries of codimension at least three.
The same difficulty one encounters in the study of manifolds 
admitting metrics of positive scalar curvature. Gromov and Lawson 
\cite{gromov-lawson80a} show that these surgeries are sufficient to
handle simply-connected manifolds (where in simply connected it is 
included that the manifold is connected).  

\begin{thm} \label{thm-bordism-GL}
\cite{gromov-lawson80a}
Suppose that $M$ is a compact simply connected spin manifold of 
dimension $\geq 5$ and suppose that $M$ is spin bordant to a 
spin manifold $N$. 
Then $M$ can be obtained from $N$ by a sequence 
of surgeries of codimension $\geq 3$.
\end{thm}

Using this theorem and the Surgery Theorem \ref{thm-spectralsurgery} 
we get the following result relating the spectrum of the Dirac operator 
for two ends of a spin bordism:

\begin{thm} \label{thm-bordism-GL-spectrum}
Let $M$ be a compact simply connected spin manifold of 
dimension $\geq 5$ and let $M$ be spin bordant to a 
spin manifold $N$.
Let $h$ be a Riemannian metric on $N$ and let $\Lambda, \varepsilon >0$ 
be such that $\pm \Lambda \not \in \spec (D_h)$. 
Then there is a Riemannian metric $g$ on $M$ such that $D_h$ and $D_g$ are 
$(\Lambda,\varepsilon )$-spectral close.   
\end{thm}

We will now show that the manifolds from 
Propositions~\ref{prop-lichnerowicz} and \ref{prop-generators} 
are sufficiently plenty to allow us to find a manifold in $\mmin$ in 
each spin bordism class. 

\begin{prop} \label{prop-mmin-bordism}
Any compact spin manifold is spin bordant to a manifold in $\mmin$. 
\end{prop}

\begin{proof}
Let $M$ be a compact spin manifold of dimension $n$ and let 
$p = \alpha(M)$. If $p \geq 0$ let $pV_n$ be the disjoint union of $p$ 
copies of the manifold $V_n$ from Proposition \ref{prop-generators} and 
if $p<0$ let $pV_n$ be the disjoint union of $-p$ copies of $-V_n$. We 
have $\alpha(M - pV_n) = \alpha(M) - p\alpha(V_n) =0$ so by Theorem B 
of Stolz \cite{stolz92a} $M - pV_n$ is spin bordant to a manifold $E$ 
which allows a positive scalar curvature metric. This means that $M$ is 
spin bordant to $E + pV_n$ which by Propositions \ref{prop-disjointunion} 
and \ref{prop-lichnerowicz} is in $\mmin$. 
\end{proof}

We are now ready to prove our main result on Conjecture B.

\begin{thm} \label{thm-simplyconnectedinmmin}
Conjecture B is true for all simply connected spin 
manifolds of dimension at least five.
\end{thm}

\begin{proof}
Suppose $M$ is a simply connected spin manifold with $\dim M \geq 5$.
From the previous proposition we know that $M$ is spin bordant to a 
manifold in $\mmin$. Theorem \ref{thm-bordism-GL} tells us that this 
bordism can be decomposed into a sequence of surgeries of codimension 
at least three. Using Proposition \ref{prop-surgery} we conclude that 
$M \in \mmin$ and that the conjecture holds for $M$. 
\end{proof}

\subsection{Non-simply connected manifolds}
In the study of manifolds allowing positive scalar curvature metrics a 
machinery has been built up to analyse when two manifolds are related by 
surgeries of codimension $\geq 3$, see for example 
\cite{rosenberg-stolz94a,rosenberg-stolz98a}. We are now going to give 
one example of how results from this area can be used to prove 
Conjecture B for certain classes of fundamental groups.

For non-simply connected manifolds one uses a refinement of the ordinary
spin bordism groups $\Omega_*^{\rm{spin}}$. Let $\pi$ be a group and let 
$B\pi$ be a classifying space for $\pi$. The bordism groups 
$\Omega_*^{\rm{spin}}(B\pi)$ consists of bordism classes of pairs $(M,f)$ 
where $M$ is a compact spin manifold and $f: M \rightarrow B\pi$ is a map. 
In particular, if $M$ has fundamental group $\pi$ we get such a pair $(M,f)$
since the universal cover of $M$ is a $\pi$-bundle over $M$ and corresponds
to a map $f : M \rightarrow B\pi$. 

In \cite{rosenberg86a} Rosenberg proves the following generalization of 
the bordism theorem \ref{thm-bordism-GL}:

\begin{thm}\label{thm-bordism-RS}
Let $M$ be a compact spin manifold of dimension $\geq 5$ with fundamental 
group $\pi$. Let $f$ be the classifying map of the universal cover of $M$. 
Suppose that $(M,f)$ is equivalent to $(N,f')$ in 
$\Omega_*^{\rm{spin}}(B\pi)$ (where $f'$ is any map $N \rightarrow B\pi$).  
Then $M$ can be obtained from $N$ by a sequence of surgeries in 
codimension $\geq 3$.
\end{thm}

For manifolds $M$ and $N$ satisfying the conditions of this theorem we 
get a bordism result for the Dirac spectrum just as Theorem 
\ref{thm-bordism-GL-spectrum}.

A finite group $\pi$ is called {\em periodic} if it has periodic cohomology,
that is if there is $d>0$ so that $H^i(\pi) = H^{i+d}(\pi)$ for all 
$i>0$. This includes cyclic groups of all orders. 

\begin{thm}\label{thm-periodic}
Suppose $M$ is a compact connected spin manifold with $\dim M \geq 5$
and fundamental group $\pi$ a periodic group of odd order. 
Then Conjecture B holds for $M$.
\end{thm}

\begin{proof}
Let $f: M \rightarrow B\pi$ be the classifying map of the 
universal cover of $M$. The class 
$(M,f) \in \Omega_*^{\rm{spin}}(B\pi)$ can be written as 
$$
(M,f) = (M,f) - (M,{\rm const}) + (M,{\rm const}).
$$ 
Composition with the constant map $B\pi \rightarrow {\rm pt}$ 
gives a map 
$\Omega_n^{\rm{spin}}(B\pi) \rightarrow \Omega_n^{\rm{spin}}({\rm pt})$
and the first two terms $(M,f) - (M,{\rm const})$ constitute an element 
in the kernel $\widetilde{\Omega}_n^{\rm{spin}}(B\pi)$ of this map. 
From Theorem 1.8 and fact (1.2) in \cite{kwasik-schultz90a} we know that 
this class can be represented by a manifold $M'$ with positive scalar 
curvature and fundamental group $\pi$. Since $M'$ has a metric of 
positive scalar curvature $M' \in \mmin$. By doing surgery on embedded 
circles we can make the third term $(M,{\rm const})$ simply connected 
without changing its class in $\Omega_*^{{\rm spin}}(B\pi)$. 
The resulting manifold $N$ is in $\mmin$ by 
Theorem \ref{thm-simplyconnectedinmmin}.

We thus have that $(M,f)$ is spin bordant to 
$(M',f') + (N, {\rm const})$ over $B\pi$. Since $\alpha (M') = 0$ 
Proposition \ref{prop-disjointunion} tells us that the disjoint union 
$M' + N \in \mmin$. Theorem \ref{thm-bordism-RS} and Proposition 
\ref{prop-surgery} finish the proof.
\end{proof}

\Remark{
The crucial point of the proof is the fact that for odd order periodic 
groups $\pi$ the kernel $\widetilde{\Omega}_n^{\rm{spin}}(B\pi)$ can 
be represented by manifolds allowing positive scalar curvature metrics.
This fact --- and the theorem --- holds for any odd order finite group 
for which the Gromov-Lawson-Rosenberg conjecture is true, see 
Section 3 of \cite{rosenberg-stolz95a}.
}

\Remark{
The torus $T^n$ with a flat metric (and the trivial spin structure)
has a large space of harmonic spinors, even though $\alpha(T^n) = 0$.
It is not clear when in general the product of two manifolds in
$\mmin$ is again in $\mmin$.
By taking products we can however prove the
following. 

{\em On the torus $T^n$ with any spin structure there is a metric for which
the kernel of the Dirac operator is trivial.
} 

Namely, from the solution of Conjecture B for manifolds of low dimensions by
Maier \cite{maier97a} we know that (for any spin structure) there is a
metric $g_3$ on $T^3$ without nontrivial harmonic spinors. Write 
$T^n = T^3 \times T^{n-3}$ and set $g_n = g_3 + g_{n-3}$, where
$g_{n-3}$ is any metric on $T^{n-3}$. 
Since $D_{g_3}$ does not have zero as an eigenvalue it follows from the 
Pythagorean formula for the eigenvalues
of the product metric $g_n$ that $D_{g_n}$ has no zero eigenvalue.
}

\providecommand{\bysame}{\leavevmode\hbox to3em{\hrulefill}\thinspace}

\end{document}